\documentclass[11pt]{article}
\usepackage[table]{xcolor}
\usepackage{graphicx}
\usepackage{amsmath,amsfonts,amssymb}
\usepackage{graphicx}
\usepackage{multirow}
\usepackage{tikz,pgf}
\usepackage{url}
\usepackage{amsmath}
\usepackage{graphicx}
\usepackage{epsfig}
\usepackage{epstopdf}
\usepackage{color}
\usepackage{enumitem}
\def\disp{\displaystyle}

\def\tto{\;{\lower 1pt \hbox{$\rightarrow$}}\kern -10pt
\hbox{\raise 2pt \hbox{$\rightarrow$}}\;}
\def\Hat{\widehat}

\def\Bar{\overline}
\def\ra{\rangle}
\def\la{\langle}
\def\ve{\varepsilon}
\def\epsilon{\varepsilon}
\def\B{\Bbb B}
\def\h{\hfill\Box}
\def\R{\Bbb R}

\def\N{\Bbb N}
\def\ox{\bar{x}}
\def\ow{\bar{w}}

\def\co{\mbox{\rm co}}

\def\dom{\mbox{\rm dom}\,}

\def\cl{\mbox{\rm cl}\,}

\def\h{\hfill\square}
\def\dn{\downarrow}
\def\ph{\varphi}
\def\emp{\emptyset}

\def\oR{\Bar{\R}}

\def\ph{\varphi}
\def\emp{\emptyset}

\def\oR{\Bar{\R}}

\setlist[enumerate,1]{itemsep=0.0ex,parsep=0.5ex,label={\rm(\alph*)},leftmargin=*,align=left}
\newcounter{lk}

\usepackage[colorlinks=true]{hyperref}

\begin{document}
\begin{center}
{\sc\bf Convex Analysis of Minimal Time and Signed Minimal Time Functions}\\[1ex]
{\sc D. V. Cuong} \footnote{Department of Mathematics, Faculty of Natural Sciences, Duy Tan University, Da Nang, Vietnam (dvcuong@duytan.edu.vn). This research is funded by Vietnam National Foundation for Science and Technology Development (NAFOSTED) under grant number 101.02-2020.20.}$^,$\footnote{Institute of Research and Development, Duy Tan University, Da Nang, Vietnam.}, {\sc B. S. Mordukhovich}\footnote{Department of Mathematics, Wayne State University, Detroit, Michigan 48202, USA (boris@math.wayne.edu). Research of this author was partly supported by the USA National Science Foundation under grants DMS-1512846 and DMS-1808978, by the USA Air Force Office of Scientific Research grant \#15RT04, and by Australian Research Council under grant DP-190100555.},
{\sc N. M. Nam}\footnote{Fariborz Maseeh Department of Mathematics and Statistics, Portland State University, Portland, OR 97207, USA (mnn3@pdx.edu). Research of this author was partly supported by the USA National Science Foundation under grant DMS-1716057.}, {\sc M. Wells}\footnote{Fariborz Maseeh Department of Mathematics and Statistics, Portland State University, Portland, OR 97207, USA (mlwells@pdx.edu).}\\[2ex]
{\bf To Alex Rubinov, in memoriam} \end{center}
\small{\bf Abstract.} In this paper we first consider the class of minimal time functions in the general setting of locally convex topological vector (LCTV) spaces. The results obtained in this framework are based on a novel notion of closedness of target sets with respect to constant dynamics. Then we introduce and investigate a new class of signed minimal time functions, which are generalizations of the signed distance functions. Subdifferential formulas for the signed minimal time and distance functions are obtained under the convexity assumptions on the given data.\\[1ex]

\newtheorem{Theorem}{Theorem}[section]
\newtheorem{Proposition}[Theorem]{Proposition}
\newtheorem{Remark}[Theorem]{Remark}
\newtheorem{Lemma}[Theorem]{Lemma}
\newtheorem{Corollary}[Theorem]{Corollary}
\newtheorem{Definition}[Theorem]{Definition}
\newtheorem{Example}[Theorem]{Example}
\renewcommand{\theequation}{\thesection.\arabic{equation}}
\normalsize\vspace*{-0.2in}

\section{Introduction}\label{intro}\vspace*{-0.1in}

The paper is dedicated to the memory of Professor Alexander Rubinov whose enormous contributions to convex and nonsmooth analysis, optimization theory, and their numerous applications are difficult to overstate; see, e.g., his book \cite{rub}.

It has been well recognized in these areas that the class of distance functions plays a crucial role in many aspects of convex and variational analysis with their applications to optimization-related and other problems; see, e.g., the books \cite{Bauschke2011,CL,m-book1,m-book,rw} with many references and discussions therein. Nevertheless, the distance function suffers from the actual drawback that it vanishes on the entire set under consideration, and thus it fails to fully characterize the geometric structure of the set near its boundary. The {\em signed distance function}, also known as the {\em oriented distance function},  was introduced to address this shortcoming. The signed distance function has been used in a broad spectrum of applications including binary classification, computer vision, shape analysis, image processing, computational physics, etc.; see the publications \cite{Wanka2013,Boczko2009,Delfour1994,H1979,ktz,H1979,Perera2015} with their bibliographies. We particularly refer the reader to \cite{Wang19} for a thorough and deep study of the class of signed distance functions on finite-dimensional Euclidean spaces with a survey of the previous achievements. It is rather surprising that the signed distance function and its generalized differential properties are much less investigated and applied in the literature in comparison with its standard distance function counterpart, even in finite dimensions.

On the other hand, recent years have witnessed great interest in the class of {\rm minimal time functions} that are obtained from the distance functions by replacing the norm function by the Minkowski gauge. This extension of the classical distance function opens up the possibility of dealing with many classes of generalized distances. Various properties of the minimal time functions have been topics of extensive research over the years including very recent developments; see, e.g., \cite{cgm,cowo,cowo1,durea1,durea2,durea3,heng,mn10,bmn10,bmn,nam-cuong,nam-zalin,nam-an2012}. Applications of the minimal time functions to variational analysis, set optimization, facility location problems, and optimal control have been largely addressed in the literature. Similarly to the distance functions, the minimal time functions vanish on the sets under consideration and thus do not capture the geometric structure of the set in question. This is the driving force for us to introduce a new class of {\em signed minimal time functions}. We define and study here this class of functions in the general setting of locally convex topological vector spaces with specifications to the cases of normed spaces and $\R^n$. Besides this goal, we also obtain some important properties of the minimal time functions defined now on locally convex topological vector spaces, in comparison with the previously known results in the case of finite-dimensional and various types of normed spaces. To proceed in this direction, we introduce in this paper a new notion of the $F$-closure of the target set with respect to the constant dynamics in the construction of the minimal time function. This notion allows us to provide a unified study of the minimal time functions in the general framework of LCTV spaces with providing simplified proofs of known results even in the normed and finite-dimensional settings. The obtained results are used in the subsequent study of the signed minimal time and distance functions.

The rest of the paper is organized as follows. In Section~\ref{sec:hor} we recall the notion of the {\em horizon cone} for convex sets in LCTV spaces and present some of its properties needed below. Section~\ref{sec:mtf} is devoted to the study of the {\em minimal time functions} associated with nonempty target subsets $\Omega\subset X$ of LCTV spaces $X$ relative to a given constant convex dynamics $F\subset X$. We define here the $F$-closure $\cl_F(\Omega)$ of a set $\Omega$ and efficiently use it to establish various properties of the minimal time functions in the general nonconvex framework. Section~\ref{sec:mtf-diff} concerns subdifferentiation of the minimal time functions generated by convex target sets in LCTV spaces. We derive here precise formulas for calculations of convex subgradient sets for such functions at given points in all the possible cases where $\ox\in\Omega$, $\ox\in\cl_F(\Omega)$, and $\ox\notin\cl_F(\Omega)$.

In Section~\ref{Signed} we introduce the class of {\em signed minimal time functions} and establish their important properties and representations in the general setting of LCTV spaces. This includes deriving useful relationships for evaluating their convex subdifferentials. Section~\ref{sec:dist} presents a new proof of the precise formula from \cite{Wang19} on calculating the convex subdifferential of the {\em signed distance function} in $\R^n$ that is based on a variational approach and subdifferential results of nonconvex variational analysis. The concluding Section~\ref{sec:concl} discuses some open questions of future research.

Throughout the paper we use the standard notation and terminology of convex and variational analysis; see, e.g., \cite{m-book1,z}. Unless otherwise stated, the space $X$ under consideration is a {\em real Hausdorff locally convex topological vector space} (LCTV space for brevity), and $X^*$ is its topological dual with the canonical pairing $\la x^*,x\ra:=x^*(x)$ for $x\in X$ and $x^*\in X^*$.\vspace*{-0.2in}

\section{The Horizon Cone}\label{sec:hor}
\setcounter{equation}{0}\vspace*{-0.1in}

This section recalls the notion of the {\em horizon cone}, which is instrumental to describe behavior of convex sets at {\em infinity}. For completeness and the reader's convenience, we present here short proofs of the statements used in what follows in LCTV spaces.

\begin{Definition}\label{horiz} Given a nonempty, closed, and convex subset $F$ of $X$ and a point $x\in F$, the {\sc horizon cone} of $F$ at $x$ is defined by
\begin{equation*}
F_\infty(x):=\big\{d\in X\;\big|\;x+td\in F\;\mbox{ for all }\;t>0\big\}.
\end{equation*}
\end{Definition}
This definition can be rewritten in the form
\begin{equation}\label{hor}
F_\infty(x)=\bigcap_{t>0}\frac{F-x}{t},
\end{equation}
which implies that $F_\infty(x)$ is a convex cone in $X$. It is also closed if $F$ is closed.

The following proposition shows that $F_\infty(x)$ does not depend on the choice of $x\in F$ provided that $F$ is a nonempty, closed, and convex subset of $X$. Thus in this case we can simply use the notation $F_\infty$ for the horizon cone of $F$.

\begin{Proposition}\label{p1} Let $F$ be a nonempty, closed, and convex subset of $X$. Then we have the equality $F_{\infty}(x_1)=F_{\infty}(x_2)$ for any $x_1,x_2\in F$.
\end{Proposition}
{\bf Proof.} It suffices to verify that $F_\infty(x_1)\subset F_\infty(x_2)$ whenever $x_1,x_2\in F$. Taking any $d\in F_\infty(x_1)$ and $t>0$, we show that $x_2+td\in F$. Consider the sequence
\begin{equation*}
x_k:=\dfrac{1}{k}\Big(x_1+ktd\Big)+\left(1-\dfrac{1}{k}\right)x_2,\quad k\in\N.
\end{equation*}
Then $x_k\in F$ for every $k$ due to $d\in F_\infty(x_1)$ and the convexity of $F$. We also have $x_k\to x_2+td$, and thus $x_2+td\in F$ since $F$ is closed. It tells us that $d\in F_\infty(x_2)$. $\h$

\begin{Corollary}\label{c1} Let $F$ be a closed and convex subset of $X$ that contains the origin. Then
\begin{equation*}
F_{\infty}=\bigcap_{t>0}tF.
\end{equation*}
\end{Corollary}
{\bf Proof.} It follows from the construction of $F_{\infty}$ and Proposition~\ref{p1} that
\begin{equation*}
F_{\infty}=F_{\infty}(0)=\bigcap_{t>0}\frac{F-0}{t}=\bigcap_{t>0}\frac{1}{t}F=\bigcap_{t>0}tF,
\end{equation*}
which readily justifies the claim. $\h$

The next proposition provides a useful sequential description of the horizon cone \eqref{hor}.

\begin{Proposition}\label{c2}
Let $F$ be a nonempty, closed, and convex subset of $X$. Then the following properties are equivalent:
\begin{enumerate}
\item[\bf(i)] We have $d\in F_{\infty}$.
\item[\bf(ii)] There exist sequences $\{t_k\}\subset(0,\infty)$ and $\{f_k\}\subset F$ such that $t_k\rightarrow 0$ and $t_kf_k\to d$.
\end{enumerate}
\end{Proposition}
{\bf Proof.} To verify implication (i)$\Longrightarrow$(ii), take $d\in F_\infty$ and fix $\bar x\in F$. It follows directly from the definition that
\begin{equation*}
\bar x+kd\in F\;\mbox{ for all }\;k\in\N.
\end{equation*}
For each $k\in\N$, this allows us to find $f_k\in F$ such that
\begin{equation*}
\bar x+kd=f_k,\;\mbox{ or equivalently, }\dfrac{1}{k}\bar x+d=\dfrac{1}{k} f_k.
\end{equation*}
Letting $t_k:=\dfrac{1}{k}$, we see that $t_kf_k\rightarrow d$ as $k\to\infty$.

To prove (ii)$\Longrightarrow$(i), suppose that there are sequences $\{t_k\}\subset(0,\infty)$ and $\{f_k\}\subset F$ with $t_k\rightarrow 0$ and $t_k f_k\to d$. Fix $x\in F$ and verify that $d\in F_\infty$ by showing that
\begin{equation*}
x+td\in F\;\mbox{ for all }\;t>0.
\end{equation*}
Indeed, for any fixed $t>0$ we have $0\le t\cdot t_k<1$ when $k$ is sufficiently large. Therefore
\begin{equation*}
(1-t\cdot t_k)x+t\cdot t_kf_k\rightarrow x+td\;\mbox{ as }\;k\to\infty.
\end{equation*}
It follows from the convexity of $F$ that every element $(1-t\cdot t_k)x+t\cdot t_kf_k$ belongs to $F$. Hence $x+td\in F$ by its closedness, and thus we get $d\in F_\infty$. $\h$\vspace*{-0.2in}

\section{Minimal Time Functions in LCTV Spaces}\label{sec:mtf}
\setcounter{equation}{0}\vspace*{-0.1in}

In this section we deal with the class of minimal time functions associated with generally nonconvex target sets in the setting of arbitrary LCTV spaces. The main goal here is to reveal important properties of such functions that do not concern their subdifferentiation. The latter topic is investigated in the next section in the convex framework.

Let $F\subset X$ be a convex set with $0\in{\rm int}(F)$, and let $\Omega$ be a nonempty subset of $X$. Define the minimal time function with the {\em target} $\Omega$ and the (constant) {\em dynamics} $F$ by
\begin{equation}\label{MTF}
\mathcal{T}_{\Omega}^F(x):=\inf\big\{t>0\;\big|\;(x+tF)\cap\Omega\ne\emp\big\},\quad x\in X.
\end{equation}
To proceed with the study of \eqref{MTF} in the general setting of LCTV spaces, we introduce the new notion of the {\em $F$-closure} of $\Omega$ by
\begin{equation}\label{f-cl}
\mbox{\rm cl}_F(\Omega):=\bigcap_{\epsilon>0}(\Omega-\epsilon F).
\end{equation}
Recall that a subset $A$ of $X$ is bounded if for any neighborhood $V$ of the origin there exists $t>0$ such that $A\subset tV$. The following proposition shows that the notion of $F$-closure reduces to the classical notion of closure if $F$ is bounded.

\begin{Proposition}\label{clFcl} Let $F\subset X$ be a convex set with $0\in\mbox{\rm int}(F)$, and let $\Omega$ be a nonempty subset of $X$. If $F$ is bounded, then $\mbox{\rm cl}_F(\Omega)=\Bar{\Omega}$.
\end{Proposition}
{\bf Proof.} Fix any $x\in\Bar{\Omega}$ and choose a neighborhood $V$ of the origin such that $V\subset F$. Then for any $\epsilon>0$ we have $x\in\Omega-\epsilon V\subset\Omega-\epsilon F$. It readily yields $x\in\mbox{\rm cl}_F(\Omega)$.

To verify the opposite inclusion, fix any $x\in\mbox{\rm cl}_F(\Omega)$ and then get $x\in\Omega-\epsilon F$ for all $\epsilon>0$. Taking any neighborhood $V$ of the origin, we have $F\subset tV$ for some $t>0$. It implies that $\epsilon F\subset V$ for $\epsilon:=1/t$, and hence $x\in\Omega-V$, i.e., $(x+V)\cap\Omega\ne\emp$. Thus $x\in\Bar{\Omega}$, which completes the proof of the proposition. $\h$

The next proposition provides a useful representation of the $F$-closure of $\Omega$ from \eqref{f-cl} via the target set $\Omega$  and the horizon cone \eqref{hor} of the dynamics $F$.

\begin{Proposition} Let $F$ be a closed and convex subset of $X$ with $0\in\mbox{\rm int}(F)$, and let $\Omega$ be a sequentially compact subset of $X$. Then we have the representation $\mbox{\rm cl}_F(\Omega)=\Omega-F_\infty$.
\end{Proposition}
{\bf Proof.} Fix any $x\in\Omega-F_\infty$ and get $x=w-d$ with some $w\in\Omega$ and $d\in F_\infty$. It follows from Corollary~\ref{c1} that $t(w-x)\in F$ for all $t>0$, and thus $x\in\Omega-\epsilon F$ for all $\epsilon>0$. This readily brings us to $x\in\mbox{\rm cl}_F(\Omega)$.

To prove the opposite inclusion, pick any $x\in\mbox{\rm cl}_F(\Omega)$ and deduce from \eqref{f-cl} that for any $k\in\N$ there exists $w_k\in\Omega$ such that $x=w_k-v_k/k$ with  $v_k\in F$. Since $\Omega$ is sequentially compact, suppose without loss of generality that the sequence $\{w_k\}$ converges to $w\in\Omega$ as $k\to\infty$. Then $v_k/k\to w-x$ as $k\to\infty$, and so Proposition~\ref{c2} ensures that $w-x\in F_\infty$. This tells us that $x\in\Omega-F_\infty$ and thus completes the proof. $\h$

Now we derive important properties of the minimal time function including its ``norm-like" description via the $F$-closure \eqref{f-cl}.

\begin{Proposition}\label{clFcl1} Let $F$ be a convex subset of $X$ with $0\in\mbox{\rm int}(F)$, and let $\Omega$ be a nonempty subset of $X$. Then for any $x\in X$ we have the following properties:
\begin{enumerate}
\item[\bf(i)] $\mathcal{T}_{\Omega}^F(x)$ is a real number.
\item[\bf(ii)] $\mathcal{T}_{\Omega}^F(x)=\inf\{t\ge 0\;|\;(x+tF)\cap\Omega\ne\emp\}$.
\item[\bf(iii)] $\mathcal{T}_{\Omega}^F(x)\ge 0$, and $\mathcal{T}_{\Omega}^F(x)=0$ if and only if $x\in\mbox{\rm cl}_F(\Omega)$.
\end{enumerate}
\end{Proposition}
{\bf Proof.} {\bf(i)} Since $0\in\mbox{\rm int}(F)$, the set $F$ is {\em absorbing}, i.e., for any $x\in X$ there exists $\nu>0$ such that $\lambda x\in F$ whenever $|\lambda|\le\nu$. Taking $w\in\Omega$, we can find $t>0$ such that $w-x\in tF$. It follows that $(x+tF)\cap\Omega\ne\emp$, and so $\mathcal{T}_{\Omega}^F(x)$ is a real number.\\[1ex]
{\bf(ii)} We obviously have the lower estimate
\begin{equation*}
\mathcal{T}_{\Omega}^F(x)\ge\gamma:=\inf\big\{t\ge 0\;\big|\;(x+tF)\cap\Omega\ne\emp\big\}.
\end{equation*}
Choose further a sequence $\{t_k\}\subset[0,\infty)$ such that $t_k\to\gamma$ as $k\to\infty$ and $(x+t_kF)\cap\Omega\ne\emp$ for every $k\in\N$. Letting $s_k:=t_k+1/k$, we get $s_k>0$ and $(x+s_kF)\cap\Omega\ne\emp$ for each $k$. Since $\{s_k\}$ converges to $\gamma$, we see that $\gamma\ge\mathcal{T}_{\Omega}^F(x)$.\\[1ex]
{\bf(iii)} It follows from the definition that $\mathcal{T}_{\Omega}^F(x)\ge 0$. Suppose that $\mathcal{T}_{\Omega}^F(x)=0$. Then for any $\epsilon>0$ there exists $0<t<\epsilon$ such that
\begin{equation*}
(x+tF)\cap\Omega\ne\emp.
\end{equation*}
Since $tF\subset\epsilon F$, we have $(x+\epsilon F)\cap\Omega\ne\emp$, which tells us that $x\in\Omega-\epsilon F$ and hence $x\in\mbox{\rm cl}_F(\Omega)$ by the construction in \eqref{f-cl}.

Conversely, supposing that $x\in\mbox{\rm cl}_F(\Omega)$ gives us $x\in\Omega-\epsilon F$, and so $(x+\epsilon F)\cap\Omega\ne\emp$ for all $\epsilon>0$. Thus $\mathcal{T}_{\Omega}^F(x)<\epsilon$, which readily yields $\mathcal{T}_{\Omega}^F(x)=0$ and completes the proof. $\h$\vspace*{0.03in}

Note that we do not assume the convexity of the target set $\Omega$ in definition \eqref{MTF} of the minimal time function. Now we show that the convexity of $\Omega$ ensures the convexity of  $\mathcal{T}_{\Omega}^F$.

\begin{Proposition}\label{mtf-conv} Let $F$ be a convex subset of $X$ with $0\in\mbox{\rm int}(F)$, and let $\Omega$ be a nonempty and convex subset of $X$. Then the minimal time function $\mathcal{T}_{\Omega}^F$ is convex on $X$.
\end{Proposition}
{\bf Proof.} Pick $x,u\in X$ and $\lambda\in(0,1)$. Then for any $\epsilon>0$ find $s>0$ and $t>0$ such that
\begin{eqnarray*}
a-x\in sF\;\mbox{ and }\;b-u\in tF\;\mbox{ for some }\;a,b\in\Omega,
\end{eqnarray*}
and that $s<\mathcal{T}_{\Omega}^F(x)+\epsilon$, $t<\mathcal{T}_{\Omega}^F(u)+\epsilon$. It implies that
\begin{eqnarray*}
\lambda a-\lambda x\in\lambda sF\;\mbox{ and }\;(1-\lambda)b-(1-\lambda)u\in(1-\lambda)t F.
\end{eqnarray*}
The convexity of the dynamics set $F$ ensures that
\begin{eqnarray*}
\big(\lambda a+(1-\lambda)b\big)-\big(\lambda x+(1-\lambda)u\big)\in\lambda s F+(1-\lambda)tF=\big(\lambda s+(1-\lambda)t\big)F.
\end{eqnarray*}
Since $\lambda a+(1-\lambda)b\in \Omega$ by the assumed convexity of the target, we get
\begin{equation*}
\mathcal{T}_{\Omega}^F\big(\lambda x+(1-\lambda)u\big)\le\lambda s+(1-\lambda)t<\lambda\mathcal{T}_{\Omega}^F(x)+(1-\lambda)\mathcal{T}_{\Omega}^F(u)+\epsilon.
\end{equation*}
Letting there $\epsilon\dn 0$ gives us the inequality
\begin{equation*}
\mathcal{T}_{\Omega}^F\big(\lambda x+(1-\lambda)u\big)\le\lambda\mathcal{T}_{\Omega}^F(x)+(1-\lambda)\mathcal{T}_{\Omega}^F(u),
\end{equation*}
which therefore justifies the convexity of the minimal time function $\mathcal{T}_{\Omega}^F$ on $X$. $\h$

To proceed further, we recall the notion of {\em Minkowski gauge function} associated with the constant dynamics set $F\subset X$ by
\begin{equation}\label{gauge}
\rho_F(x):=\inf\big\{t>0\;\big|\;x\in tF\big\},\quad x\in X.
\end{equation}
Note that in the case where $X$ is a normed space and $F$ is the closed unit ball $\B\subset X$, the Minkowski gauge \eqref{gauge} reduces to the classical distance function $d(x;F)$ associated with $F$. The next proposition gives us a representation of the minimal time function \eqref{MTF} by using the Minkowski gauge in the general LCTV setting.

\begin{Proposition}\label{prop5} Let $F$ be a convex subset of $X$ with $0\in\mbox{\rm int}(F)$, and let $\Omega$ be a nonempty subset of $X$. Then we have the representation
\begin{equation}\label{mtf-min}
\mathcal{T}_{\Omega}^F(x)=\inf_{w\in\Omega}\rho_F(w-x),\quad x\in X.
\end{equation}
\end{Proposition}
{\bf Proof.} Fix any $t>0$ such that $(x+tF)\cap\Omega\ne\emp$ and then have $a-x\in tF$ for some $a\in\Omega$. It follows from \eqref{gauge} that $\rho_F(a-x)\le t$, and thus $\inf_{w\in \Omega}\rho_F(w-x)\le t$, which tells us that $\inf_{w\in\Omega}\rho_F(w-x)\leq \mathcal{T}_{\Omega}^F(x)$. To verify the opposite inequality, fix any $t>0$ with $\inf_{w\in \Omega}\rho_F(w-x)<t$. Then there exists $a\in\Omega$ satisfying $\rho_F(a-x)<t$. This yields $(a-x)/t\in F$, and so $a-x\in tF$. Then we get $\mathcal{T}_{\Omega}^F(x)\leq t$ and hence arrive at the estimate $\mathcal{T}_{\Omega}^F(x)\le\inf_{w\in\Omega}\rho_F(w-x)$, which completes the proof.
$\h$

As a consequence of Proposition~\ref{prop5}, we show now that the minimal time functions with constant dynamics $F$ and respective target sets $\Omega$ and $\overline{\Omega}$ agree.

\begin{Proposition}\label{prop6} Let the assumptions of Proposition \ref{prop5} hold. Then we have
\begin{equation*}
\mathcal{T}_{\Omega}^F(x)=\mathcal{T}_{\overline{\Omega}}^F(x)\;\mbox{ for all }\;x\in X.
\end{equation*}
\end{Proposition}
{\bf Proof.}  By Proposition \ref{prop5}, it suffices to show that
\begin{equation*}
\disp\inf_{w\in\Omega}\rho_F(w-x)=\inf_{w\in\overline{\Omega}}\rho_F(w-x).
\end{equation*}
It is obvious that $\inf_{w\in\overline{\Omega}}\rho_F(w-x)\le\inf_{w\in\Omega}\rho_F(w-x)$. To prove the opposite inequality, pick any $t>0$ with $t>\inf_{w\in\overline{\Omega}}\rho(w-x)$. Then there exists $a\in\overline{\Omega}$ such that $\rho_F(a-x)<t$, which implies that $a-x\in tF$. Since $F$ contains a neighborhood of the origin, for any $\epsilon>0$ we get that $(a+\epsilon F)\cap\Omega\ne \emp$. Thus there exist $f\in F$ and $b\in\Omega$ with $a+\epsilon f=b$. It follows from the convexity of $F$ that $b-x\in tF+\epsilon F=(t+\epsilon)F$. Thus
\begin{equation*}
\rho_F(b-x)\le t+\epsilon\;\mbox{ and }\;\inf_{w\in \Omega}\rho_F(w-x)\le t+\epsilon.
\end{equation*}
Remembering that $\epsilon>0$ was chosen arbitrarily, we obtain  $\inf_{w\in\Omega}\rho_F(w-x)\le t$ and therefore arrive at the estimate
\begin{equation*}
\disp\inf_{w\in\Omega}\rho_F(w-x)\le\disp\inf_{w\in\overline{\Omega}}\rho_F(w-x),
\end{equation*}
which completes the proof of the proposition. $\h$

Having in hand the obtained representation from Proposition~\ref{prop5} allows us to derive some properties of the minimal time function \eqref{MTF} from those for the Minkowski gauge \eqref{gauge}.
In particular, the following lemma for \eqref{gauge} is used below.

\begin{Lemma}\label{mlinear1} If $F$ is a convex subset of $X$ with $0\in\mbox{\rm int}(F)$, then the Minkowski gauge function $\rho_F$ is continuous on $X$.
\end{Lemma}
{\bf Proof.} Fixing any $\epsilon>0$, we deduce from definition \eqref{gauge} that
\begin{equation*}
\epsilon F\subset\rho_F^{-1}([0,\epsilon]).
\end{equation*}
For any open set $V$ of $\R$ containing $0$, find $\epsilon>0$ such that $[0,\epsilon]\subset V$. Then
\begin{equation*}
\epsilon F\subset\rho_F^{-1}(V).
\end{equation*}
Since $\epsilon F$ is a neighborhood of the origin, the function $\rho_F$ is continuous at $0\in X$. Then the subadditivity of $\rho_F$  yields its continuity on the whole space $X$. $\h$

Define further the {\em polar} of $F$ by
\begin{equation*}
F^{\circ}:=\big\{ x^*\in X^*\;\big|\;\la x^*,x\ra\le 1\;\mbox{ for all }\;x\in F\big\}.
\end{equation*}
If $X$ is a normed space, we also define the {\em norm }of $F^\circ$ by
\begin{equation*}
\|F^{\circ}\|:=\mbox{\rm sup}\big\{\|x^*\|\;\big|\;x^*\in F^{\circ}\big\}.
\end{equation*}

The above constructions lead us to the next proposition.

\begin{Proposition}\label{MTcont} Let $F$ be a convex set in $X$ with $0\in\mbox{\rm int}(F)$, and let $\Omega$ be a nonempty subset of $X$. Then we have the estimate
\begin{equation}\label{mintime1}
\big|\mathcal{T}_{\Omega}^F(x)-\mathcal{T}_{\Omega}^F(y)\big|\le\max\big\{\rho_F(y-x),\rho_F(x-y)\big\}\;\mbox{ for all }\;x,y\in X,
\end{equation}
which ensures, in particular, that the minimal time function $\mathcal{T}_{\Omega}^F$ is continuous on $X$. If in addition $X$ is a normed space, then $\mathcal{T}_{\Omega}^F$ is Lipschitz continuous on $X$ with constant $\|F^\circ\|$.
\end{Proposition}
{\bf Proof.} Fix any $x,y\in X$ and any $q\in\Omega$. Then Proposition~\ref{prop5} and the subadditivity of the Minkowski gauge $\rho_F$ imply that
\begin{eqnarray*}
\begin{array}{ll}
\mathcal{T}_{\Omega}^F(x)&=\disp\inf_{w\in\Omega}\rho_F(w-x)\\
&\le\rho_F(q-x)=\rho_F(q-y+y-x)\\
&\leq \rho_F(q-y)+\rho_F(y-x).
\end{array}
\end{eqnarray*}
Taking the infimum with respect to $q\in\Omega$ and interchanging the roles of $x$ and $y$ verify \eqref{mintime1}. Fix further any $x_0\in X$ and consider the function
\begin{equation*}
\widehat{\rho}_F(x):=\max\big\{\rho_F(x-x_0),\rho_F(x_0-x)\big\}\;\mbox{ for }\;x\in X.
\end{equation*}
It follows from Lemma~\ref{mlinear1} that $\widehat{\rho}_F$ is continuous on $X$. Thus whenever $\epsilon>0$ there exists a neighborhood $V$ of $x_0$ such that
\begin{equation*}
\big|\mathcal{T}_{\Omega}^F(x)-\mathcal{T}_{\Omega}^F(x_0)\big|\le\widehat{\rho}_F(x)=\big|\widehat{\rho}_F(x)-\widehat{\rho}_F(x_0)\big|<\epsilon\;\mbox{ for all }\;x\in V,
\end{equation*}
which justifies the continuity of $\mathcal{T}_{\Omega}^F$ around $x_0$ and therefore on the entire space $X$.

To verify the last statement of the proposition in the case where $X$ is a normed space, observe that $\|F^\circ\|<\infty$ and $\rho_F(x)\le\|F^\circ\|\cdot\|x\|$ on $X$.
Then it follows from \eqref{mintime1} that
\begin{equation*}
\left|\mathcal{T}_{\Omega}^F(x)-\mathcal{T}_{\Omega}^F(y)\right|\le\|F^\circ\|\cdot\|x-y\|\;\mbox{ for all }\;x,y\in X,
\end{equation*}
which shows that $\mathcal{T}_{\Omega}^F$ is Lipschitz continuous with constant $\|F^\circ\|$. $\h$\vspace*{-0.2in}

\section{Subgradients of Minimal Time Functions in LCTV Spaces}\label{sec:mtf-diff}
\setcounter{equation}{0}\vspace*{-0.1in}

As we know from Proposition~\ref{mtf-conv}, the minimal time function \eqref{MTF} is convex provided that both target and dynamics sets are convex in the general setting of LCTV spaces $X$. This section is devoted to convex subdifferentiation of $\mathcal{T}_{\Omega}^F$ at all the possible situations where $\ox\in\Omega$, $\ox\in\cl_F(\Omega)$, and $\ox\notin\cl_F(\Omega)$. When the underlying space $X$ is Banach (or merely normed), some of the obtained results can be found in \cite{bmn10} and \cite{nam-an2012}.

Recall that, given a nonempty convex set $\Omega$ in an LCTV space $X$, the {\em normal cone} to $\Omega$ at $\ox\in\Omega$ is defined by the formula
\begin{equation}\label{nc}
N(\ox;\Omega):=\big\{x^*\in X\;\big|\;\la x^*,x-\ox\ra\le 0\;\text{ for all }\;x\in\Omega\big\}
\end{equation}
with $N(\ox;\Omega):=\emp$ if $\ox\notin\Omega$. Considering further an extended-real-valued convex function $\ph\colon X\to\oR:=(-\infty,\infty]$, the {\em subdifferential} (collections of subgradients) of $\ph$ at $\ox\in\dom(\ph):=\{x\in X\;|\;\ph(x)<\infty\}$ is given by
\begin{equation}\label{sub}
\partial\ph(\ox):=\big\{x^*\in X^*\;\big|\;\la x^*,x-\ox\ra\le\varphi(x)-\varphi(\ox)\;\mbox{ for all }\;x\in X\big\}.
\end{equation}
Recall also the {\em support function} $\sigma_F\colon X^*\to\oR$ of $F$ that is defined by
\begin{equation}\label{supp}
\sigma_F(x^*):=\sup\big\{\la x^*,f\ra\;\big|\;f\in F\big\},\quad x^*\in X^*.
\end{equation}

Note that Proposition~\ref{clFcl1}(i) tells us that the domain of the minimal time function \eqref{MTF} is the whole space $X$ if $0\in{\rm int}(F)$. The first result of this section establishes a precise formula for calculating the subdifferential of $\mathcal{T}_{\Omega}^F$ at the target points $\ox\in\Omega$ via the normal cone to $\Omega$ and the support function of the constant dynamics $F$.

\begin{Theorem}\label{inset1} Let $F$ be a convex set in $X$ with $0\in\mbox{\rm int}(F)$, and let $\Omega$ a nonempty convex subset of $X$. Then for any $\bar x\in\Omega$ we have
\begin{equation}\label{sub1}
\partial\mathcal{T}^F_{\Omega}(\ox)=N(\ox;\Omega)\cap C^*,
\end{equation}
where the set $C^*\subset X^*$ is defined via the support function \eqref{supp} of $F$ by
\begin{equation}\label{c*}
C^*:=\big\{x^*\in X^*\;\big|\;\sigma_F(-x^*)\le 1\big\}.
\end{equation}
Furthermore, we get the equalities
\begin{equation}\label{sub1a}
\partial\mathcal{T}^F_{\Omega}(\ox)=N\big(\ox;\mbox{\rm cl}_F(\Omega)\big)\cap C^*=N(\ox;\Omega)\cap C^*.
\end{equation}
\end{Theorem}
{\bf Proof.} Fix any $x^*\in\partial\mathcal{T}^F_ {\Omega}(\ox)$ and write by definition \eqref{sub} that
\begin{equation}\label{convexity}
\la x^*,x-\bar x\ra\le\mathcal{T}^F_{\Omega}(x)-\mathcal{T}^F_{\Omega}(\ox)\;\mbox{ for all }\;x\in X.
\end{equation}
Since $\mathcal{T}^F_{\Omega}(x)=0$ whenever $x\in\Omega$, it follows that
\begin{equation*}
\la x^*,x-\ox\ra\le 0\;\mbox{ for all }\;x\in\Omega,
\end{equation*}
which implies by the normal cone definition \eqref{nc} that $x^*\in N(\ox;\Omega)$. Fixing now any $f\in F$ and $t>0$, we deduce from (\ref{convexity}) that
\begin{equation*}
\la x^*,(\ox-tf)-\ox\ra\le\mathcal{T}^F_{\Omega}(\ox-tf)\le t,
\end{equation*}
where the last inequality holds due to $((\ox-tf)+tF)\cap\Omega\ne\emp$. This ensures therefore that
\begin{equation*}
\la x^*,-f\ra\le 1\;\mbox{ for all }\;f\in F,
\end{equation*}
and so $x^*\in C^*$ by the construction of $C^*$ and definition \eqref{supp} of the support function. Thus we arrive at the inclusion $\partial\mathcal{T}^F_{\Omega}(\ox)\subset N(\ox;\Omega)\cap C^*$.

To verify the opposite inclusion in \eqref{sub1}, pick any $x^*\in N(\ox;\Omega)\cap C^*$ and then get $\la x^*,-f\ra\le 1$ for all $f\in F$ together with
\begin{equation*}
\la x^*,x-\ox\ra\le 0\;\mbox{ for all }\;x\in\Omega.
\end{equation*}
Fix further any $u\in X$ and for every $\ve>0$ find $t>0$, $f\in F$, and $\omega\in\Omega$ such that
\begin{equation*}
\mathcal{T}^F_{\Omega}(u)\le t<\mathcal{T}^F_{\Omega}(u)+\ve\;\mbox{ and }\;u+tf=\omega.
\end{equation*}
This readily gives us the relationships
\begin{eqnarray*}
\begin{array}{ll}
\la x^*,u-\ox\ra&=\la x^*,u-\omega\ra+\la x^*,\omega-\ox\ra\\
&\le\la x^*,-tf\ra\le t<\mathcal{T}^F_{\Omega}(u)+\ve=\mathcal{T}^F_{\Omega}(u)-\mathcal{T}^F_{\Omega}(\ox)+\ve,
\end{array}
\end{eqnarray*}
which yield \eqref{convexity} since $\ve$ was chosen arbitrarily, and thus justify formula \eqref{sub1}.

Let us finally prove the fulfillment of \eqref{sub1a} that reduces by \eqref{sub1} to verifying the second equality therein. We obviously have the inclusion $N\big(\ox;\mbox{\rm cl}_F(\Omega)\big)\cap C^*\subset N(\ox;\Omega)\cap C^*$ due to $\Omega\subset\mbox{\rm cl}_F(\Omega)$. To verify the opposite one, pick any $x^*\in N(\ox;\Omega)\cap C^*$ and $x\in\mbox{\rm cl}_F(\Omega)$. Then we get from \eqref{f-cl} that $x\in\Omega-tF$ for all $t>0$. Representing $x=\omega_t-td_t$ with some $\omega_t\in\Omega$ and $d_t\in F$ for each $t>0$ tells us that
\begin{equation*}
\la x^*,x-\ox\ra =\la x^*,\omega_t-td_t-\ox\ra=\la x^*,\omega_t-\ox\ra+t\la x^*,-d_t\ra\le t,
\end{equation*}
which ensures by passing to the limit as $t\dn 0$ that $x^*\in N(\ox;\mbox{\rm cl}_F(\Omega))\cap C^*$ and thus completes the proof of the theorem. $\h$\vspace*{0.03in}

The next theorem derives a similar subdifferential formula for $\mathcal{T}_{\Omega}^F$ at points $\ox\in\mbox{\rm cl}_F(\Omega)$.

\begin{Theorem}\label{outset1} Let $\ox\in\mbox{\rm cl}_F(\Omega)$ at the setting of Theorem~{\rm\ref{inset1}}. Then we have
\begin{equation}\label{sub2}
\partial\mathcal{T}^F_{\Omega}(\ox)=N\big(\ox;\mbox{\rm cl}_F(\Omega)\big)\cap C^*,
\end{equation}
where the set $C^*\subset X^*$ is taken from \eqref{c*}.
\end{Theorem}
{\bf Proof.} Picking any $x^*\in\partial\mathcal{T}^F_{\Omega}(\ox)$, deduce from \eqref{convexity} and the equality $\mathcal{T}^F_{\Omega}(x)=0$ for all $x\in\mbox{\rm cl}_F(\Omega)$ that
$x^*\in N(\ox;\mbox{\rm cl}_F(\Omega))$. Since $\mbox{\rm cl}_F(\Omega)\subset\Omega-F$, we represent $\ox$ as $\ow-\bar{f}$ with $\ow\in\Omega$ and $\bar{f}\in F$. Fixing now any $f\in F$ and $t>0$, denote $x:=\ow-tf$ and obtain similarly to the proof of Theorem~\ref{inset1} that
\begin{equation*}
\la x^*,(\ow-tf)-(\ow-\bar{f})\ra=\la x^*,-tf+\bar{f}\ra\le\mathcal{T}_{\Omega}^F(\ow-tf)\le t,
\end{equation*}
which clearly ensures that
\begin{equation*}
\la x^*,-f+\bar{f}/t\ra\le 1\;\mbox{ for all }\;t>0.
\end{equation*}
Letting $t\to\infty$ and using \eqref{c*} give us $x^*\in C^*$ and hence verify the inclusion
\begin{equation*}
\partial\mathcal{T}^F_{\Omega}(\ox)\subset N\big(\ox;\mbox{\rm cl}_F(\Omega)\big)\cap C^*.
\end{equation*}

To justify the opposite inclusion in \eqref{sub2}, pick any $x^*\in N\big(\ox;\mbox{\rm cl}_F(\Omega)\big)\cap C^*$. Taking into account Proposition~\ref{clFcl1}(i), for every
$u\in X$ find $t\in[0,\infty)$, $\omega\in\Omega$, and $f\in F$ with
\begin{equation*}
\mathcal{T}^F_{\Omega}(u)\le t<\mathcal{T}^F_{\Omega}(u)+\ve\;\mbox{ and }\;u+tf=\omega.
\end{equation*}
Since $\Omega\subset\mbox{\rm cl}_F(\Omega)$, we deduce from the definitions that
\begin{equation*}
\la x^*,u-\ox\ra=\la x^*,\omega-\ox\ra+t\la x^*,-f\ra\le t<\mathcal{T}^F_{\Omega}(u)+\ve=\mathcal{T}^F_{\Omega}(u)-\mathcal{T}^F_{\Omega}(\ox)+\ve.
\end{equation*}
This therefore ensures that $x^*\in\partial\mathcal{T}^F_{\Omega}(\ox)$ and thus verifies the inclusion ``$\supset$" in \eqref{sub2}. In this way we complete the proof of the theorem. $\h$\vspace*{0.03in}

In the rest of this section we consider the most challenging case for subdifferentiation of the minimal time function \eqref{MTF} at the points  $\ox\notin\mbox{\rm cl}_F(\Omega)$. Two precise and different formulas for computing the subdifferential $\partial\mathcal{T}^F_{\Omega}(\ox)$ at such points are derived below. The first subdifferential formula involves the {\em expansions} of $\Omega$ defined by
\begin{equation}\label{exp}
\Omega_r:=\big\{x\in X\;\big|\;\mathcal{T}^F_{\Omega}(x)\le r\big\}\;\mbox{ for any }\;r>0.
\end{equation}

To proceed, we begin with the following two technical lemmas, where we do not assume that $0\in{\rm int}(F)$ and hence do not guarantee that ${\rm dom}(\mathcal{T}^F_{\Omega})=X$.

\begin{Lemma}\label{distanceestimate} Let $F$ be a convex set in $X$, let $\Omega$ be a nonempty subset of $X$, and let $r>0$. Suppose $x\notin\Omega_r$ and $\mathcal{T}^F_{\Omega}(x)<\infty$. Then we have
\begin{equation*}
\mathcal{T}^F_{\Omega}(x)=\mathcal{T}^F_{\Omega_r}(x)+r.
\end{equation*}
\end{Lemma}
{\bf Proof.} Fix any $t>0$ such that $(x+tF)\cap\Omega_r\ne\emp$ and then find $f_1\in F$ and $u\in\Omega_r$ with $x+tf_1=u$. It follows from $u\in\Omega_r$ that $\mathcal{T}^F_{\Omega}(u)\le r$. This tells us by \eqref{f-cl}  that for any $\ve>0$ there exists $s>0$ such that $s<r+\ve$ and $(u+s F)\cap\Omega\ne\emp$. Consequently we get $\omega\in\Omega$ and $f_2\in F$ such that $u+sf_2=\omega$. Therefore
\begin{equation*}
\omega=u+sf_2=(x+tf_1)+sf_2\in x+(t+s)F
\end{equation*}
by the convexity of the constant dynamics $F$. This clearly implies that
\begin{equation*}
\mathcal{T}^F_{\Omega}(x)\le t+s\le t+r+\ve.
\end{equation*}
Since $\ve>0$ was chosen arbitrarily, we deduce that $\mathcal{T}^F_{\Omega}(x)\le t+r$, and hence
\begin{equation*}
\mathcal{T}^F_{\Omega}(x)\le\mathcal{T}^F_{\Omega_r}+r.
\end{equation*}

To verify the opposite inequality, denote $\gamma:=\mathcal{T}^F_{\Omega}(x)$ and get $r<\gamma$ since $x\notin\Omega_r$.
For any $\ve>0$ we find $t\in[0,\infty)$, $f\in F$, and $\omega\in\Omega$ with
\begin{equation*}
\gamma\le t<\gamma+\ve\;\mbox{ and }\;x+tf=\omega.
\end{equation*}
This gives us the relationships
\begin{equation*}
\omega=x+tf=x+(t-r)f+rf\in x+(t-r)f+rF
\end{equation*}
and therefore verifies the estimate
\begin{equation*}
\mathcal{T}^F_{\Omega}(x+(t-r)f)\le r.
\end{equation*}
The latter implies that $x+(t-r)f\in\Omega_r$ and also that $x+(t-r)f\in x+(t-r)F$. Hence we have $\mathcal{T}^F_{\Omega_r}(x)\le t-r\le\gamma-r+\ve$. Since $\ve>0$ was chosen arbitrarily, this brings us to
\begin{equation*}
r+\mathcal{T}^F_{\Omega_r}\le\gamma=\mathcal{T}^F_{\Omega}(x)
\end{equation*}
and thus completes the proof of the lemma. $\h$\vspace*{0.03in}

The second lemma provides an estimate of the shifted minimal time function.

\begin{Lemma}\label{lm3} Let $F$ be a convex set in $X$, let $\Omega$ be a nonempty subset of $X$, and let $x\in\mbox{\rm dom}(\mathcal{T}^F_{\Omega})$. Then for $t\ge 0$ and any $f\in F$ we have
\begin{equation*}
\mathcal{T}^F_{\Omega}(x-tf)\le\mathcal{T}^F_{\Omega}(x)+t.
\end{equation*}
\end{Lemma}
{\bf Proof.} Pick any $\ve>0$ and find $s\ge 0$ such that
\begin{equation*}
\mathcal{T}^F_{\Omega}(x)\le s<\mathcal{T}^F_{\Omega}(x)+\ve\;\mbox{ and }\;(x+sF)\cap\Omega\ne\emp.
\end{equation*}
Then $(x-tf+tF+sF)\cap\Omega\ne\emp$ and hence $(x-tf+(t+s)F)\cap\Omega\ne\emp$ telling us that
\begin{equation*}
\mathcal{T}^F_{\Omega}(x-tf)\le t+s\le\mathcal{T}^F_{\Omega}(x)+t+\ve.
\end{equation*}
Passing there to the limit as $\ve\dn 0$ verifies the claim. $\h$\vspace*{0.03in}

Now we are ready to obtain the first formula to calculate the subdifferential of $\mathcal{T}^F_{\Omega}$ at points $\ox$ outside of $\mbox{\rm cl}_F(\Omega)$. The obtained result invokes the normal cone \eqref{nc} to the expansion \eqref{exp} of $\Omega$ with $r=\mathcal{T}^F_{\Omega}(\ox)$ at $\ox\in\Omega_r$.

\begin{Theorem}\label{outset} In the setting of Theorem~{\rm\ref{inset1}}, suppose that $\bar x\notin\mbox{\rm cl}_F(\Omega)$, and let $r:=\mathcal{T}^F_{\Omega}(\ox)>0$. Then we have the subdifferential formula
\begin{equation}\label{sub3}
\partial\mathcal{T}^F_{\Omega}(\ox)=N(\ox;\Omega_r)\cap S^*\;\mbox{ with }\;S^*:=\big\{x^*\in X^*\;\big|\;\sigma_F(-x^*)=1\big\}.
\end{equation}
\end{Theorem}
{\bf Proof.} Pick any $x^*\in\partial\mathcal{T}^F_{\Omega}(\ox)$ and conclude similarly to the proof of Theorem~\ref{inset1} that $\sigma_F(-x^*)\le 1$ and $x^*\in N(\ox;\Omega_r)$. Let us now show that $\sigma_F(-x^*)=1$. Having in mind that
\begin{equation}\label{convexity1}
\la x^*,x-\ox\ra\le\mathcal{T}^F_{\Omega}(x)-\mathcal{T}^F_{\Omega}(\ox)\;\mbox{ for all }\;x\in X,
\end{equation}
fix any $\ve\in(0,r)$ and find $t\in\R$, $f\in F$, and $\omega\in\Omega$ satisfying
\begin{equation*}
r\le t<r+\ve^2\;\mbox{ and }\;\omega=\ox+tf.
\end{equation*}
We can write $\omega=\ox+\ve f+(t-\ve)f$ and thus get the estimate $\mathcal{T}^F_{\Omega}(\ox+\ve f)\le t-\ve$. Applying \eqref{convexity1} with $x=\ox+\ve f$ gives us the inequalities
\begin{equation*}
\la x^*,\ve f\ra\le\mathcal{T}^F_{\Omega}(\ox+\ve f)-\mathcal{T}^F_{\Omega}(\ox)\le t-\ve-r\le\ve^2-\ve,
\end{equation*}
which ensure the fulfillment of the estimates
\begin{equation*}
1-\ve\le\la-x^*,f\ra\le\sigma_F(-x^*),
\end{equation*}
They imply by passing to the limit as $\ve\dn 0$ that $\sigma_F(-x^*)\ge 1$, and thus $x^*\in S^*$. This verifies the inclusion $\partial\mathcal{T}^F_{\Omega}(\ox)\subset N(\ox;\Omega_r)\cap S^*$.

To justify the opposite inclusion in \eqref{sub3}, take any $x^*\in N(\ox;\Omega_r)$ for which $\sigma_F(-x^*)=1$ and then show that the subgradient inequality \eqref{convexity1} is satisfied. It follows from Theorem~\ref{inset1} that $x^*\in\partial\mathcal{T}^F_{\Omega_r}(\ox)$, and thus
\begin{equation*}
\la x^*,x-\ox\ra\le\mathcal{T}^F_{\Omega_r}(x)\;\mbox{ for all }\;x\in X.
\end{equation*}
Fix any $x\in X$ and first consider the case where $t:=\mathcal{T}^F_{\Omega}(x)>r$. Then it follows from Lemma~\ref{distanceestimate} that the desired condition \eqref{convexity1} holds. In the remaining case where $t\le r$, for any $\ve>0$ we choose $f\in F$ such that $\la x^*,-f\ra>1-\ve$. Employing Lemma~\ref{lm3} ensures that $\mathcal{T}^F_{\Omega}(x-(r-t)f)\le r$, and hence $x-(r-t)f\in\Omega_r$. Since $x^*\in N(\ox;\Omega_r)$, we get
\begin{equation*}
\la x^*,x-(r-t)f-\ox\ra\le 0,
\end{equation*}
which clearly yields the relationships
\begin{equation*}
\la x^*,x-\ox\ra\le\la x^*,f\ra(r-t)\le(1-\ve)(t-r)=(1-\ve)\big(\mathcal{T}^F_{\Omega}(x)-\mathcal{T}^F_{\Omega}(\ox)\big).
\end{equation*}
Since $\ve>0$ was chosen arbitrarily, we arrive at \eqref{convexity1} and thus complete the proof. $\h$\vspace*{0.03in}

The last result of this section establishes another representation of the subdifferential of $\mathcal{T}^F_{\Omega}$ at points $\bar x\notin\mbox{\rm cl}_F(\Omega)$ that uses the subdifferential of the Minkowski gauge \eqref{gauge} of the dynamics $F$ and the normal cone to the target $\Omega$ at the generalized projection points instead of the target expansion $\Omega_r$. Now we impose more assumptions on $\Omega$ that ensure, in particular, that the $F$-closure $\mbox{\rm cl}_F(\Omega)$ reduces to the standard closure $\Bar\Omega$ by Proposition~\ref{clFcl}. Also the class of LCTV spaces should be specified. Recall that a (Hausdorff) LCTV space is {\em semi-reflexive} if the canonical map into its bidual is surjective; see, e.g., \cite{mv}.

\begin{Theorem}\label{subofs} In addition to the assumptions of Theorem~{\rm\ref{inset1}}, suppose that $X$ is semi-reflexive and $\Omega$ is closed and bounded. Then for any $\ox\in X$ we have the representation
\begin{equation}\label{sub4}
\partial\mathcal{T}^F_{\Omega}(\ox)=\big(-\partial\rho_F(\bar\omega-\ox)\big)\cap N(\bar\omega;\Omega)
\end{equation}
for any $\bar\omega\in\Pi_F(\ox;\Omega)$, where the nonempty generalized projection is defined by
\begin{equation}\label{proj}
\Pi_F(\ox;\Omega):=\big\{\omega\in\Omega\;\big|\;\mathcal{T}^F_{\Omega}(\ox)=\rho_F(\omega-\ox)\big\}.
\end{equation}
\end{Theorem}
{\bf Proof.} Let us first observe that the generalized projection \eqref{proj} is a nonempty set. Indeed, Lemma~\ref{mlinear1} tells us that the Minkowski gauge $\rho_F$ is a continuous function on $X$ under the imposed convexity of $F$ with $0\in{\rm int}(F)$. The convexity of $\rho_F$ ensures that $\rho_F$ is weakly lower semicontinuous on $X$. Due to the semi-reflexivity of $X$ and the assumptions imposed on $\Omega$, we conclude that $\Omega$ is weakly compact in $X$; see, e.g., \cite[Proposition~23.18]{mv}. Thus the infimum in \eqref{mtf-min} is realized by the Weierstrass existence theorem in the weak topology of $X$. This justifies the nonemptiness of $\Pi_F(\ox;\Omega)$. Now representation \eqref{sub4} readily follows from the proof of \cite[Theorem~7.3]{bmn10}, and thus we are done. $\h$\vspace*{-0.2in}

\section{Signed Minimal Time Functions}\label{Signed}
\setcounter{equation}{0}\vspace*{-0.1in}

In this section we introduce and study a new class of functions that seems to be important in convex and variational analysis. Without loss of generality in the study of the following class of function, assume that the target sets $\Omega$ are {\em proper}, i.e., they are nonempty together with their complements $\Omega^c$. Given a target set $\Omega\subset X$ and a constant dynamics $F\subset X$ in an LCTV space $X$, the {\em signed minimal time function} is defined by
\begin{equation}\label{smt}
\Delta_{\Omega}^F(x):=\begin{cases}
\mathcal{T}_{\Omega}^F(x)&\mbox { if }\;x\notin\Omega,\\
-\mathcal{T}_{\Omega^c}^F(x)&\mbox{ if }\;x\in\Omega
\end{cases}
\end{equation}
via the minimal time functions \eqref{MTF} of $\Omega$ and $\Omega^c$. Consider also the associated extended-real-valued function $\mu_{\Omega}^F\colon X\to\oR$ given by
\begin{equation}\label{mu}
\mu_{\Omega}^F(x):=\begin{cases}
\infty&\mbox{ if }\;x\in\Omega^c,\\
-\mathcal{T}_{\Omega^c}^F(x)&\mbox{ if }\;x\in\Omega.
\end{cases}
\end{equation}

First we establish some properties of function \eqref{mu}.

\begin{Proposition}\label{Mu} Let $F$ be a convex set with $0\in\mbox{\rm int}(F)$, and let $\Omega,\Omega_1,\Omega_2$ be proper subsets of $X$. Then we have the following properties of \eqref{mu}:
\begin{enumerate}
\item[\bf(i)] If $\Omega_1\subset\Omega_2$, then $\mu_{\Omega_2}^F(x)\le\mu_{\Omega_1}^F(x)$ for all $x\in X$.
\item[\bf(ii)] $\big\{x\in X\;\big|\;\mu_{\Omega}^F(x)\le 0\big\}=\Omega$.
\item[\bf(iii)] If $F$ is bounded, then $\{x\in X\;|\;\mu_{\Omega}^F(x)<0\}=\mbox{\rm int}(\Omega)$.
\item[\bf(iv)] If $\Omega$ is closed, then $\mu_{\Omega}^F$ is lower semicontinuous on $X$.
\end{enumerate}
\end{Proposition}
\textbf{Proof.} Properties (i) and (ii) follow directly from the definition. To verify (iii), it suffices to show that $\mathcal{T}_{\Omega^c}^F(x)>0$ if and only if $x\in\mbox{\rm int}(\Omega)$.
If $\mathcal{T}_{\Omega^c}^F(x)>0$, we obviously have $x\in\Omega$. The only fact to be checked is that $x\notin\mbox{\rm bd}(\Omega)$. Arguing by contraposition, suppose that $x\in\mbox{\rm bd}(\Omega)$. Since $F$ contains a neighborhood of the origin, for any $k\in\N$ we get
\begin{equation*}
\Big(x+\frac{1}{k}F\Big)\cap\Omega^c\ne\emp.
\end{equation*}
It follows that $\mathcal{T}_{\Omega^c}^F(x)\le 1/k$ for all $k\in\mathbb{N}$, and hence $\mathcal{T}_{\Omega^c}^F(x)=0$, a contradiction.

Conversely, suppose that $x\in\mbox{\rm int}(\Omega)$. To verify that $\mu_{\Omega}^F(x)<0$, assume the contrary., i.e., $\mathcal{T}_{\Omega^c}^F(x)=0$. Using the imposed boundedness of $F$ and employing Propositions~\ref{clFcl} and \ref{clFcl1}(iii) tell us that $x\in\overline{\Omega^c}$, which is a contradiction.

Finally, let us verify (iv) by taking any $\alpha\in\R$ and considering the level set
\begin{equation*}
\mathcal{L}_{\alpha}:=\big\{x\in X\;\big|\;\mu_{\Omega}^F(x)\le\alpha\big\}=\big\{x\in\Omega\;\big|\;\mathcal{T}_{\Omega^c}^F(x)\ge\alpha\big\},
\end{equation*}
which is closed in $X$ by the closedness of $\Omega$ and the continuity of $\mathcal{T}_{\Omega^c}^F$ due to Proposition~\ref{MTcont}. This ensures therefore that $\mu_{\Omega}^F$ is lower semicontinuous on $X$. $\h$\vspace*{0.03in}

The next proposition shows that the convexity of function $\mu_{\Omega}^F$ from \eqref{mu} is equivalent to the convexity of the set $\Omega$ therein.

\begin{Proposition}\label{mu-convex} Let $F$ and $\Omega$ satisfy the assumptions of Proposition~{\rm\ref{Mu}}. Then $\mu_{\Omega}^F$ is convex if and only if $\Omega$ is convex.
\end{Proposition}
\textbf{Proof.} If $\mu_{\Omega}^F$ is convex, then we immediately get from Proposition~\ref{Mu}(ii) that $\Omega$ is a convex set. Conversely, assume that $\Omega$ is convex and observe from definition \eqref{mu} that the convexity of $\mu_{\Omega}^F$ reduces to the fact that
\begin{equation}\label{mu-conv}
\mathcal{T}_{\Omega^c}^F\big(\lambda x_1+(1-\lambda)x_2\big)\ge\lambda \mathcal{T}_{\Omega^c}^F(x_1)+(1-\lambda)\mathcal{T}_{\Omega^c}^F(x_2)
\end{equation}
for all $\lambda\in(0,1)$ and $x_1,x_2\in\Omega$. To verify \eqref{mu-conv}, suppose the contrary and then find $\lambda\in(0,1)$ and $x_1,x_2\in\Omega$ such that
\begin{equation*}
\mathcal{T}_{\Omega^c}^F\big(\lambda x_1+(1-\lambda)x_2\big)<\lambda\mathcal{T}_{\Omega^c}^F(x_1)+(1-\lambda)\mathcal{T}_{\Omega^c}^F(x_2).
\end{equation*}
It gives us $0<t<\lambda\mathcal{T}_{\Omega^c}^F(x_1)+(1-\lambda)\mathcal{T}_{\Omega^c}^F(x_2)$ for which
\begin{equation*}
\big[\big(\lambda x_1+(1-\lambda)x_2\big)+tF\big]\cap\Omega^c\ne\emp.
\end{equation*}
Thus we find $w\in F$ satisfying the inclusion
\begin{equation}\label{mu-conv1}
\big(\lambda x_1+(1-\lambda)x_2\big)+tw\in\Omega^c.
\end{equation}
Consider further the points
\begin{equation*}
u_1:=x_1+\frac{tw}{\gamma}\mathcal{T}_{\Omega^c}^F(x_1)\;\mbox{ and }\;u_2:=x_2+\frac{tw}{\gamma}\mathcal{T}_{\Omega^c}^F(x_2),
\end{equation*}
where $\gamma:=\lambda\mathcal{T}_{\Omega^c}^F(x_1)+(1-\lambda)\mathcal{T}_{\Omega^c}^F(x_2)$. Taking into account that $t<\gamma$, let us show now that $u_1,u_2\in\Omega$.
Indeed, for $\mathcal{T}_{\Omega^c}^F(x_1)=0$ we readily have $u_1=x_1\in\Omega$. In the case where $\mathcal{T}_{\Omega^c}^F(x_1)>0$, suppose on the contrary that
$u_1\notin\Omega$ and hence get $u_1\in\Omega^c$. Since $\big(x_1+\frac{t}{\gamma}\mathcal{T}_{\Omega^c}^F(x_1)F\big)\cap\Omega^c\ne\emp$, we obtain that
\begin{equation*}
\mathcal{T}_{\Omega^c}^F(x_1)\le\frac{t}{\gamma}\mathcal{T}_{\Omega^c}^F(x_1)<\mathcal{T}_{\Omega^c}^F(x_1),
\end{equation*}
a clear contradiction, which shows therefore that $u_1\in\Omega$. Similarly we have $u_2\in\Omega$. Using now the convexity of $\Omega$ tells us that
\begin{equation*}
\lambda u_1+(1-\lambda)u_2=\lambda x_1+(1-\lambda)x_2+tw\in\Omega,
\end{equation*}
which contradicts \eqref{mu-conv1} and thus completes the proof of the proposition. $\h$\vspace*{0.03in}

Now we are ready to derive some properties of the signed minimal time function \eqref{smt}. The following proposition shows that---in contrast to the minimal time function---its signed counterpart allows us to fully describe the interior, closure, boundary, and exterior of $\Omega$.

\begin{Proposition}\label{smt-desc} Let $F$ be a bounded, convex subset of $X$ with $0\in\mbox{\rm int}(F)$, and let $\Omega\subset X$ be proper. Then we have the following descriptions:
\begin{enumerate}
\item[\bf(i)] $\mbox{\rm int}(\Omega)=\big\{x\in X\;\big|\;\Delta_{\Omega}^F(x)<0\big\}$.
\item[\bf(ii)] $\overline{\Omega}=\big\{x\in X\;\big|\;\Delta_{\Omega}^F(x)\le 0\big\}$.
\item[\bf(iii)] $\mbox{\rm bd}(\Omega)=\big\{x\in X\;\big|\;\Delta_{\Omega}^F(x)=0\big\}$.
\item[\bf(iv)] $\mbox{\rm int}(\Omega^c)=\big\{x\in X\;\big|\;\Delta_{\Omega}^F(x)>0\big\}$.
\end{enumerate}
\end{Proposition}
\textbf{Proof.} Assertions (i) and (iv) follow from Proposition~\ref{Mu} and the definitions in \eqref{smt} and \eqref{mu}, while assertion (iii) is a consequence of (i) and (iv).
Finally, we get (ii) as a consequence of (i) and (iii). $\h$\vspace*{0.03in}

Our next result verifies the continuity of $\Delta_{\Omega}^F$ in the general LCTV setting and its Lipschitz continuity in normed spaces.

\begin{Proposition}\label{prop3} Let $F$ be a convex subset of $X$ with  $0\in\mbox{\rm int}(F)$, and let $\Omega$ be a proper subset of $X$. Then the signed minimal time function $\Delta_{\Omega}^F$ is continuous on $X$. If in addition $X$ is a normed space and if $F$ is bounded, then the function $\Delta_{\Omega}^F$ is Lipschitz continuous on $X$ with the uniform Lipschitz constant $\|F^{\circ}\|$.
\end{Proposition}
\textbf{Proof.} It is easy to observe from definition \eqref{smt} that
\begin{equation}\label{decom}
\Delta_{\Omega}^F(x)=\mathcal{T}_{\Omega}^F-\mathcal{T}_{\Omega^c}^F\;\mbox{ for all }\;x\in X
\end{equation}
since the minimal time function \eqref{MTF} vanishes at points in its target set. This tells us by Proposition~\ref{MTcont} that $\Delta_{\Omega}^F$ is a continuous function on $X$. Considering further the case where $X$ is a normed space, pick any $x,u\in X$ and deduce from Proposition~\ref{MTcont} that both functions $\mathcal{T}_{\Omega}^F$ and $\mathcal{T}_{\Omega^c}^F$ are Lipschitz continuous on $X$ with the uniform Lipschitz constant $\|F^{\circ}\|$. To verify the Lipschitz continuity of $\Delta_{\Omega}^F$ on $X$ with the same constant $\|F^{\circ}\|$, taking any points $x,u\in X$ and taking \eqref{decom} into account, it suffices to consider the case where $x\in\Omega$ and $u\in\Omega^c$. Then we get the equalities
\begin{align*}
|\Delta_{\Omega}^F(x)-\Delta_{\Omega}^F(u)|&=|\mathcal{T}_{\Omega}^F(x)-\mathcal{T}_{\Omega^c}^F(x)-\mathcal{T}_{\Omega}^F(u)+\mathcal{T}_{\Omega^c}^F(u)|\\
&=\mathcal{T}_{\Omega^c}^F(x)+\mathcal{T}_{\Omega}^F(u).
\end{align*}
Define now the continuous function $\psi\colon[0,1]\to\R$ by
\begin{equation*}
\psi(t):=\mathcal{T}_{\Omega}^F\big(tx+(1-t)u\big)-\mathcal{T}_{\Omega^c}^F\big(tx+(1-t)u\big),\quad t\in[0,1].
\end{equation*}
Since $\psi(0)\cdot\psi(1)\le 0$, the intermediate point theorem gives us $z\in(x,u)$ such that $\mathcal{T}_{\Omega^c}(z)=\mathcal{T}_{\Omega}^F(z)$.  It is easy to see that in order for this equation to hold, both $\mathcal{T}_{\Omega}^F(z)$ and $\mathcal{T}_{\Omega^c}^F(z)$ must be equal to $0$. Proposition~\ref{smt-desc} and \eqref{decom} imply then that $z\in\mbox{\rm bd}(\Omega)$. Since $z\in\overline{\Omega}\cap\overline{\Omega^c}$, Propositions \ref{prop5} and \ref{prop6} together show that
\begin{equation*}
\mathcal{T}_{\Omega^c}^F(x)+\mathcal{T}_{\Omega}^F(u)\le\rho_F(z-x)+\rho_F(z-u)\le\|F^{\circ}\|\cdot\|x-z\|+\|F^{\circ}\|\cdot\|z-u\|=\|F^{\circ}\|\cdot\|x-u\|,
\end{equation*}
which verifies the Lipschitz continuity of $\Delta_{\Omega}^F$ on $X$ with Lipschitz constant $\|F^{\circ}\|$. $\h$\vspace*{0.03in}

To obtain the next result, recall that {\em infimal convolution} of two functions $f,g\colon X\to\oR$ is
\begin{equation*}
(f\oplus g)(x):=\inf\big\{f(y)+g(x-y)\;\big|\;y\in X\big\},\quad x\in X.
\end{equation*}

The following theorem expresses the signed minimal time function \eqref{smt} as the infimal convolution of the function $\mu_{\Omega}^F$ from \eqref{mu} and the Minkowski gauge \eqref{gauge}.

\begin{Theorem}\label{convol} Let $F$ be a convex subset of $X$ with $0\in\mbox{\rm int}(F)$, and let $\Omega$ be a proper subset of $X$. Then the signed minimal time function \eqref{smt} is represented as
\begin{equation}\label{convol-eq}
\Delta_{\Omega}^F(x)=(\mu_{\Omega}^F\oplus\rho_F)(x)\;\mbox{ for all }\;x\in\Omega.
\end{equation}
If we assume in addition that the set $F$ is symmetric, then \eqref{convol-eq} holds for all $x\in X$.
\end{Theorem}
\textbf{Proof.} Let us first consider the case where $x\in\Omega$. Then we have
\begin{align*}
(\mu_{\Omega}^F\oplus\rho_F)(x)&=\inf\big\{\mu_{\Omega}^F(y)+\rho_F(x-y)\;\big|\;y\in X\big\}\\
&=\inf\big\{\mu_{\Omega}^F(y)+\rho_F(x-y)\;\big|\;y\in\Omega\big\}\\
&=\inf\big\{-\mathcal{T}_{\Omega^c}^F(y)+\rho_F(x-y)\;\big|\;y\in\Omega\big\}\\
&\le-\mathcal{T}_{\Omega^c}^F(x)+\rho_F(x-x)=-\mathcal{T}_{\Omega^c}^F(x).
\end{align*}
Taking any $y\in\Omega$, deduce from the proof of Proposition~\ref{MTcont} that
\begin{equation*}
-\mathcal{T}_{\Omega^c}^F(x)\le-\mathcal{T}_{\Omega^c}^F(y)+\rho_F(x-y).
\end{equation*}
It follows therefore that
\begin{equation*}
-\mathcal{T}_{\Omega^c}^F(x)\le\inf\big\{\mu_{\Omega}^F(y)+\rho_F(x-y)\;\big|\;y\in\Omega\big\}=(\mu_{\Omega}^F\oplus\rho_F)(x),
\end{equation*}
which verifies representation \eqref{convol-eq} on $\Omega$.

To prove the last statement of the theorem, suppose that $F$ is symmetric and pick $x\in\Omega^c$. Then for any $y\in\Omega$ we have the estimate
\begin{equation*}
-\mathcal{T}_{\Omega^c}^F(y)+\rho_F(y-x)\le\rho_F(y-x),
\end{equation*}
which readily implies that
\begin{equation*}
(\mu_{\Omega}^F\oplus\rho_F)(x)\le\inf\big\{\rho_F(y-x)\;\big|\;y\in\Omega\big\}=\mathcal{T}_{\Omega}^F(x).
\end{equation*}
Picking now $y\in\Omega$ and following the proof of Proposition~\ref{prop3}, we find $z\in\mbox{\rm bd}(\Omega)$ such that
\begin{equation*}
\rho_F(x-z)+\rho_F(z-y)=\rho_F(x-y).
\end{equation*}
This brings us to the relationships
\begin{equation*}
\mathcal{T}_{\Omega}^F(x)+\mathcal{T}_{\Omega^c}^F(y)\le\rho_F(z-x)+\rho_F(z-y)=\rho_F(x-y),
\end{equation*}
which give us the estimate
\begin{equation*}
\mathcal{T}_{\Omega}^F(x)\le-\mathcal{T}_{\Omega^c}^F(y)+\rho_F(x-y)\;\mbox{ for all }\;y\in\Omega.
\end{equation*}
Thus we arrive at $\mathcal{T}_{\Omega}^F(x)\le(\mu_{\Omega}^F\oplus\rho_F)(x)$ and complete the proof of the theorem. $\h$\vspace*{0.03in}

It is not hard to present an example showing that representation \eqref{convol-eq} may not hold for all $x\in\Omega$ if the set $F$ is not symmetric.

The following corollary is a direct consequence of Theorem~\ref{convol}.

\begin{Corollary}\label{smt-conv} Let the target set $\Omega$ be convex in addition to all the assumptions of Theorem~{\rm\ref{convol}}. Then the signed minimal time function $\Delta_{\Omega}^F$ is convex on $X$.
\end{Corollary}
\textbf{Proof.} Proposition~\ref{mu-convex} tells us that the function $\mu_{\Omega}^F$ from \eqref{mu} is convex under the imposed assumptions. Furthermore, it is well known that the infimal convolution of convex functions is convex. Thus the claimed convexity of $\Delta_{\Omega}^F$ on $X$ follows from Theorem~\ref{convol}. $\h$\vspace{0.03in}

The concluding part of this section concerns convex subdifferentiation of the signed minimal time functions \eqref{smt} in the general setting of LCTV spaces. We start with subdifferentiation of the auxiliary function \eqref{mu}, which is convex under our assumptions. Although the next propositions deals only with $\mu_{\Omega}^F$, its proof employs some properties of the signed minimal time function $\Delta_{\Omega}^F$ established above.

\begin{Proposition}\label{sub-mu} Let $F$ be a convex, symmetric subset of $X$ with $0\in\mbox{\rm int}(F)$, and let the target set $\Omega\subset X$ be proper, convex, and closed in $X$.
Then following hold:
\begin{enumerate}
\item[\bf(i)] $\partial\mu_{\Omega}^F(\ox)\ne\emp$ for all $\ox\in\Omega$.
\item[\bf(ii)] If $\ox\in\mbox{\rm bd}(\Omega)$, then we have the inclusion
\begin{equation*}
\partial\mu_{\Omega}^F(\ox)\subset N(\ox;\Omega).
\end{equation*}
\end{enumerate}
\end{Proposition}
\textbf{Proof.} To verify (i), it suffices to consider the case where $\ox\in\mbox{\rm bd}(\Omega)$. Observe that the closedness of $\Omega$ guarantees that $\mbox{\rm bd}(\Omega)\subset \Omega=\mbox{\rm dom}(\mu^F_\Omega)$. Along with $\mu_{\Omega}^F$, consider the corresponding signed minimal time function $\Delta_{\Omega}^F$, which is convex and continuous on $X$ due to Corollary~\ref{smt-conv} and Proposition~\ref{prop3}, respectively. Basic convex analysis tells us that $\partial\Delta_{\Omega}^F(\ox)\ne\emp$. Pick any $x^*\in\partial\Delta_{\Omega}^F(\ox)$ and get by \eqref{sub} and Proposition~\ref{smt-desc}(iii) that
\begin{equation*}
\la x^*,x-\ox\ra\le\Delta_{\Omega}^F(x)-\Delta_{\Omega}^F(\ox)=\Delta_{\Omega}^F(x)\;\mbox{ for all }\;x\in X.
\end{equation*}
This clearly implies the conditions
\begin{equation*}
\la x^*,x-\ox\ra\le\Delta_{\Omega}^F(x)=-\mathcal{T}_{\Omega^c}^F(x)=\mu_{\Omega}^F(x)-\mu_{\Omega}^F(\ox)\;\mbox{ whenever }\;x\in\Omega,
\end{equation*}
and hence $x^*\in\partial\mu_{\Omega}^F(\ox)$, which justifies assertion (i).

To verify (ii), take any $x^*\in\partial\mu_{\Omega}^F(\ox)$ and deduce from the definitions that
\begin{equation*}
\la x^*,x-\ox\ra\le\mu_{\Omega}^F(x)-\mu_{\Omega}^F(\ox)=-\mathcal{T}_{\Omega^c}^F(x)\le 0\;\mbox{ for all }\;x\in\Omega.
\end{equation*}
This gives us $x^*\in N(\ox;\Omega)$ and thus completes the proof of the proposition. $\h$\vspace*{0.03in}

The last result here provides a precise formula for representing the subdifferential of $\Delta_{\Omega}^F$ via that of $\mu_{\Omega}^F$ at boundary points of $\Omega$.

\begin{Proposition}\label{sub-smt} Let all the assumptions of Proposition~{\rm\ref{sub-mu}} be satisfied. Then we have
\begin{equation*}
\partial\Delta_{\Omega}^F(\ox)=\partial\mu_{\Omega}^F(\ox)\cap F^\circ\;\mbox{ for every }\;\ox\in\mbox{\rm bd}(\Omega).
\end{equation*}
\end{Proposition}
\textbf{Proof.} Fix any $\ox\in{\rm bd}(\Omega)$. Theorem~\ref{convol} gives us the representation $\Delta_{\Omega}^F(x)=(\mu_{\Omega}^F\oplus\rho_F)(x)$ for all $x\in X$, where the infimal convolution is {\em exact} at the reference point $\ox$ in the sense that $\Delta_{\Omega}^F(\ox)=\mu_{\Omega}^F(\ox)+\rho_F(0)$, since all the terms therein are zero. Applying the subdifferential rule for infimal convolutions of convex functions defined on LCTV spaces (see, e.g., \cite[Corollary~2.4.7]{z}), we get the equalities
\begin{equation*}
\partial\Delta_{\Omega}^F(\ox)=\partial(\mu_{\Omega}^F\oplus\rho_F)(\ox)=\partial\mu_{\Omega}^F(\ox)\cap\partial\rho_F(0)=\partial\mu_{\Omega}^F(\ox)\cap F^\circ
\end{equation*}
and thus verify the claimed subdifferential representation for \eqref{smt}. $\h$

It follows from the combination of Proposition~\ref{sub-mu} and Proposition~\ref{sub-smt} that, under the assumptions imposed therein, the subdifferential of the signed minimal time function \eqref{smt} on an LCTV space $X$ admits the upper estimate
\begin{equation*}
\partial\Delta_{\Omega}^F(\ox)\subset N(\ox;\Omega)\cap F^\circ\;\mbox{ for every }\;\ox\in\mbox{\rm bd}(\Omega).
\end{equation*}\vspace*{-0.35in}

\section{Subgradients of Signed Distance Functions}\label{sec:dist}
\setcounter{equation}{0}\vspace*{-0.1in}

Now we are going to provide a more detailed subdifferential study of the particular case of \eqref{smt} corresponding to $F=\B\subset\R^n$ with the Euclidean norm on $\R^n$. In this case, the signed minimal time function reduces to the {\em signed distance function} defined by
\begin{equation}\label{s-dist}
\Hat d(x;\Omega):=\begin{cases}
d(x;\Omega)&\mbox{if }\;x\in\Omega^c,\\
-d(x;\Omega^c)&\mbox{if }\;x\in\Omega,
\end{cases}
\end{equation}
where $d(x;\Theta)$ stands for the standard {\em distance function} associated with a set $\Theta\subset\R^n$ by
\begin{equation}\label{dist}
d(x;\Theta):=\inf\big\{\|w-x\|\;\big|\;w\in\Theta\big\},\quad x\in\R^n.
\end{equation}
Generalized differential properties of the distance function \eqref{dist} have been well investigated in the cases of convex and nonconvex sets in finite-dimensional and various infinite-dimensional settings; see, e.g., \cite{CL,m-book1,m-book,bmn,rw} and the references therein. This has been also largely done for the class of {\em marginal/optimal value functions} given in the form
\begin{equation}\label{marg}
\vartheta(x):=\inf\big\{\ph(x,y)\;\big|\;y\in G(x)\big\},\quad x\in X,
\end{equation}
which plays a highly important role in many aspects of convex and variational analysis, optimization, and their numerous applications; see the aforementioned books. However, it is not the case of the signed distance function \eqref{s-dist}, which cannot be reduced to the marginal function form \eqref{marg}.

The main goal of this section is to derive a precise (equality type) formula for calculating the convex subdifferential of the signed distance function at each point of the Euclidean space $\R^n$. Such a result was established in the recent paper \cite[Theorem~3.8]{Wang19} by using an involved machinery of convex analysis. Here we provide a different proof, which reduces the calculation of the subdifferential \eqref{sub} of the {\em convex} signed distance function \eqref{s-dist} to the calculation of the limiting subdifferential of the {\em nonconvex} standard distance function \eqref{dist} by involving tools of variational analysis.

To proceed, for each $x\in\Omega$ consider the set
\begin{equation}\label{Q}
Q_{\Omega}(x):=\begin{cases}
\Pi(x;\overline{\Omega})&\mbox{ if }\;x\in\Omega^c,\\
\Pi(x;\overline{\Omega^c})&\mbox{ if }\;x\in\Omega,
\end{cases}
\end{equation}
where $\Pi(x;\Theta)$ stands for the Euclidean projection of $x$ to $\Theta$ defined by
\begin{equation}\label{pr}
\Pi(x;\Theta):=\big\{w\in\Theta\;\big|\;\|x-w\|=d(x;\Theta)\big\}.
\end{equation}
It is well known that $\Pi(x;\Theta)$ is nonempty if $\Theta$ is closed while being a singleton when $\Theta$ is closed and convex. In the case of $F=\B$ under consideration, the function $\mu_F^\Omega$ from \eqref{mu} reduces to $\mathcal{\theta}\colon\R^n\to\oR$ given by
\begin{equation}\label{th}
\mathcal{\theta}(x):=\begin{cases}
\infty&\mbox{ if }\;x\in\Omega^c,\\
-d(x;\Omega^c)&\mbox{ if }\;x\in\Omega.
\end{cases}
\end{equation}

Before deriving the main result of this section, let us present two useful lemmas. The first one is taken from \cite[Corollary~2.5.3]{BW}.

\begin{Lemma}\label{limiting}
Let $f\colon\R^n\to\R$ be a convex function, let $S$ be any subset of Lebesgue measure zero in $\R^n$, and let $D_f$ be the set of points in $\R^n$ at which $f$ fails to be differentiable. Then
\begin{equation*}
\partial f(\ox)={\rm co}\big\{\lim\nabla f(x_k)\;\big|\;x_k\to\ox,\;x_k\notin S\cup D_f\big\}\;\mbox{ for all }\;\ox\in\R^n,
\end{equation*}
where the symbol `${\rm co}$' stands for the convex hull of a set.
\end{Lemma}

For the second lemma we provide a simple proof.

\begin{Lemma}\label{reversenormal} Let $\Omega$ be a proper, closed, and convex subset of $\R^n$. If $\ow\in Q_\Omega(\ox)$ as defined in \eqref{Q} with $\ox\in\mbox{\rm int}(\Omega)$, then we have the inclusion $\ow-\ox\in N(\ow;\Omega)$.
\end{Lemma}
{\bf Proof.} Fix any $x\in\overline{\Omega^c}$ and write for the Euclidean norm that
\begin{equation*}
\|x-\ox\|^2=\|x-\ow\|^2-2\la\ox-\ow,x-\ow\ra+\|\ox-\ow\|^2.
\end{equation*}
Since $\ow\in Q_{\Omega}(\ox)=\Pi(\ox;\overline{\Omega^c})$, we get $\|\ox-x\|^2-\|\ox-\ow\|^2\ge 0$, which implies that
\begin{equation*}
\la\ox-\ow,x-\ow\ra\le\dfrac{1}{2}\|x-\ow\|^2\;\mbox{ for all }\;x\in\overline{\Omega^c}.
\end{equation*}
Now let us pick any $x\in\mbox{\rm int}(\Omega)$ and verify that $\ow+t(\ow-x)\notin\Omega$ whenever $t>0$. Indeed, supposing on the contrary that $\ow+t(\ow-x)=\omega\in\Omega$ yields
\begin{equation*}
\ow=\frac{t}{t+1}x+\frac{1}{t+1}\omega\in\mbox{\rm int}(\Omega),
\end{equation*}
a contradiction. Thus we arrive at the estimate
\begin{eqnarray*}
\begin{array}{ll}
t\la\ox-\ow,\ow-x\ra&=\big\la\ox-\ow,\big(\ow+t(\ow-x)\big)-\ow\big\ra\\
&\le\dfrac{1}{2}\big\|\big(\ow+t(\ow-x)\big)-\ow\big\|^2=\dfrac{1}{2}t^2\big\|\ow-x\big\|^2\;\mbox{ for all }\;t>0.
\end{array}
\end{eqnarray*}
Letting there $t\dn 0$ gives us $\la\ow-\ox,x-\ow\ra\le 0$ and thus verifies that $\ow-\ox\in N(\ow;\Omega)$. $\h$\vspace*{0.03in}

To proceed further, we need to recall the two subdifferential constructions of variational analysis dealing with {\em locally Lipschitzian} functions on $\R^n$.

Let $f\colon\R^n\to\oR$ be locally Lipschitzian around $\ox\in\dom(f)$. Then the (Mordukhovich) {\em limiting subdifferential} of $f$ at $\ox$ is defined by
\begin{equation}\label{lim-sub}
\partial_M f(\ox):=\disp\Big\{v\in\R^n\;\Big|\;\exists\,x_k\to\ox,\;v_k\to v,\;\liminf_{x\to x_k}\frac{f(x)-f(x_k)-\la v_k,x-x_k\ra}{\|x-x_k\|}\le 0\Big\}.
\end{equation}
We refer the reader to the books \cite{m-book1,m-book,rw} and the bibliographies therein for systematic studies  and applications of this construction in finite and infinite dimensions.

The (Clarke) {\em generalized gradient} of a locally Lipschitzian function $f$ at a given point $\ox$ is defined and comprehensively studied in \cite{CL} in the Banach space setting, while we employ here its equivalent representation in $\R^n$ via the convex hull of the limiting subdifferential \eqref{lim-sub}; see, e.g., \cite[Theorem~3.57]{m-book1}:
\begin{equation}\label{cl}
\partial_C f(\ox)={\rm co}\,\partial_M f(\ox).
\end{equation}
If $f$ is convex, then both subdifferentials \eqref{lim-sub} and \eqref{cl} reduce to the subdifferential of convex analysis \eqref{sub}.

Recall \cite[Proposition~2.3.1]{CL} the plus-minus symmetry property
\begin{equation}\label{sym}
\partial_C(-f)(\ox)=-\partial_C f(\ox)
\end{equation}
of the generalized gradient \eqref{cl} used in what follows.

Now we are ready to derive the following precise calculation formula for the convex subdifferential of the signed distance functions \eqref{s-dist} at each point $\ox\in\R^n$ with the different representations at $\ox\in\Omega^c$, $\ox\in\mbox{\rm bd}(\Omega)$, and $\ox\in\mbox{\rm int}(\Omega)$. Its proof given below exploits the subdifferential constructions \eqref{lim-sub} and \eqref{cl} for nonconvex Lipschitzian functions. In particular, we are going to use the limiting subdifferential calculation
\begin{equation}\label{sub-dist}
\partial_M d(\ox;\Theta):=\disp\bigg\{\frac{\ox-\Pi(\ox;\Theta)}{d(\ox;\Theta)}\bigg\},\quad\ox\notin\Theta,
\end{equation}
via the Euclidean projection \eqref{pr} of $\ox$ to a closed set $\Theta$; see \cite[Example~8.53]{rw} and \cite[Theorem~1.33]{m-book} for this formula and its different proofs.

\begin{Theorem}\label{sdist-sub} Let $\Omega$ be a proper, closed, and convex subset of $\R^n$. Then we have
\begin{equation}\label{SDF}
\partial\Hat d(\ox;\Omega)=\begin{cases}
\bigg\{\dfrac{\ox-Q_{\Omega}(\ox)}{\Hat d(\ox;\Omega)}\bigg\}&\mbox{ if }\;\ox\in\Omega^c,\\
\bigg\{\dfrac{\ox-\textnormal{co}\,Q_{\Omega}(\ox)}{\Hat d(\ox;\Omega)}\bigg\}&\mbox{ if }\;\ox\in\mbox{\rm int}(\Omega),\\
\textnormal{co}\big(\mathbb{\Bbb S}\cap N(\ox;\Omega)\big)&\mbox{ if }\;\ox\in\mbox{\rm bd}(\Omega),
\end{cases}
\end{equation}
where the set $Q_\Omega(\ox)$ is taken from \eqref{Q}, and where $\mathbb{\Bbb S}$ stands for the unit sphere of $\R^n$.
\end{Theorem}
{\bf Proof.} Let us split the proof into the three steps corresponding to the position of the reference point $\ox\in\R^n$ with respect to the set $\Omega$ in \eqref{SDF}.

{\bf Step~1:} $\ox\in\Omega^c$. In this case we have that $\Hat d(x;\Omega)=d(x;\Omega)$ for all $x\in\Omega^c$, where $\Omega^c$ is an open set in $\R^n$.  Using the well-known result of convex analysis on subdifferentiation of the distance function to a convex set (see, e.g., \cite[Theorem~2.39]{bmn}) tells us that
\begin{equation*}
\partial\Hat d(\ox;\Omega)=\partial d(\ox;\Omega)=\bigg\{\dfrac{\ox-\Pi(\ox;\Omega)}{d(\ox;\Omega)}\bigg\}=\bigg\{\dfrac{\ox-Q_{\Omega}(\ox)}{\Hat d(\ox;\Omega)}\bigg\},\quad\ox\in\Omega^c,
\end{equation*}
which verifies the subdifferential representation of \eqref{SDF} in this case.

{\bf Step~2:} $\ox\in\mbox{\rm int}(\Omega)$. In this case we deduce from the definitions that
\begin{equation*}
\Hat d(x;\Omega)=-d(x;\Omega^c)=-d(x;\overline{\Omega^c})=\theta(x)\;\mbox{ for all }\;x\in\mbox{\rm int}(\Omega),
\end{equation*}
and thus the function $\theta$ from \eqref{th} is convex. Then applying the equalities in \eqref{cl} and \eqref{sym} together with the subdifferential calculations \eqref{sub-dist} for the distance function gives us the following relationships:
\begin{eqnarray*}
\begin{array}{ll}
\partial\Hat d(\ox;\Omega)&=\partial\theta(\ox)=\partial_C\theta(\ox)\\
&=-\partial_C d(x;\overline{\Omega^c})=-\mbox{\rm co}\big(\partial_M d(\ox;\overline{\Omega^c})\big)\\
&=-\mbox{\rm co}\bigg(\dfrac{\ox-Q(\ox)}{d(\ox;\overline{\Omega^c})}\bigg)=\bigg\{\dfrac{\ox-{\rm co}\,Q_{\Omega}(\ox)}{\Hat d(\ox;\Omega)}\bigg\}.
\end{array}
\end{eqnarray*}
Thus we arrive at the claimed formula \eqref{SDF} in the case where $\ox\in\mbox{\rm int}(\Omega)$.

{\bf Step~3:} $\ox\in\mbox{\rm bd}(\Omega)$. Since the signed distance function $\Hat d(\cdot;\Omega)$ is convex and Lipschitz continuous on $\R^n$, and since the boundary set $\mbox{\rm bd}(\Omega)$ has Lebesgue measure zero on $\R^n$, we have by Lemma~\ref{limiting} that
\begin{equation}\label{sdc}
\partial\Hat d(\ox;\Omega)=\co\big\{\lim\nabla f(x_k)\;\big|\;x_k\rightarrow\ox,\;x_k\notin\mbox{\rm bd}(\Omega)\cup D_f\big\},
\end{equation}
where we denote $f(x):=\Hat d(x;\Omega)$ on $\R^n$ for convenience.

Select a sequence $\{x_k\}\subset\R^n$ such that $x_k\notin\mbox{\rm bd}(\Omega)\cup D_f$ for every $k\in\N$, that $x_k\to\ox$ as $k\to\infty$, and that the limit $v:=\lim_{k\to\infty}\nabla f(x_k)$ exists. Consider first the case where for some $k_0\in\N$ we have that $x_k\notin\Omega$ whenever $k\ge k_0$. It clearly follows in this case that
\begin{equation*}
\nabla f(x_k)=\dfrac{x_k-\Pi(x_k)}{d(x_k;\Omega)}\in N(\Pi(x_k);\Omega)\cap\Bbb S\;\mbox{ for all }\;k\ge k_0,
\end{equation*}
and hence we get $v\in N(\ox;\Omega)\cap\Bbb S$ by passing to the limit as $k\to\infty$ in view of \eqref{nc}.

It remains to consider the case where there is a subsequence of $\{x_k\}$ such that (without relabeling) $x_k\in\mbox{\rm int}(\Omega)$ as $k\in\N$. Thus we are in the situation of Step~2, which gives us
\begin{equation*}
\nabla f(x_k)=\dfrac{x_k-w_k}{\Hat d(x_k;\Omega)}\;\mbox{ with }\;w_k:={\rm co}\,Q_{\Omega}(x_k),
\end{equation*}
where the set ${\rm co}\,Q_{\Omega}(x_k)$ is a singleton, and where $\|\nabla f(x_k)\|=1$ for each $k\in\N$. The latter ensures therefore that $v\in\Bbb S$. Furthermore, by Lemma~\ref{reversenormal} we have
\begin{equation*}
\dfrac{x_k-{\rm co}\,Q_{\Omega}(x_k)}{\Hat d(x_k;\Omega)}=\dfrac{w_k-x_k}{d(x_k;\Omega^c)}\in N(w_k;\Omega).
\end{equation*}
Since $w_k\to\ox$ as $k\to\infty$, it follows from the above constructions that $v\in N(\ox;\Omega)$, and hence $v\in N(\ox;\Omega)\cap\Bbb S$. Now we get by \eqref{sdc} that $\partial\Hat d(\ox;\Omega)\subset\co(N(\ox;\Omega)\cap\Bbb S)$.

To verify the opposite inclusion, pick any $v\in N(\ox;\Omega)\cap\Bbb S$ and let $x_k:=\ox+v/k$ for all $k\in\N$. It is not hard to check that $\nabla f(x_k)=v$ for all $k$, and hence $v\in\partial\Hat d(\ox;\Omega)$. This justifies formula \eqref{SDF} for $\ox\in{\rm bd}(\Omega)$ and thus completes the proof of the theorem. $\h$\vspace*{0.03in}

We conclude this section with the construction of an example showing that the finite dimension of the Euclidean space $X=\R^n$ is essential for the fulfillment of Theorem~\ref{sdist-sub}. In fact, the following example demonstrates that the subdifferential formula \eqref{SDF} fails in any infinite-dimensional separable Hilbert space.

\begin{Example}\label{exa} {\rm Let $X$ be an infinite-dimensional separable Hilbert space with an orthonormal basis $\{e_k\}_{k\in\N}$. Consider the set
\begin{equation*}
\Omega:=\big\{x\in X\;\big|\;\la x,e_k\ra\le k^{-1}\;\mbox{ for all }\;k\in\N\big\},
\end{equation*}
which is clearly closed and convex. We are going to show that $\ox=0\in\Omega$ is a boundary point of the set $\Omega$, and that
\begin{equation*}
\partial\Hat d(0;\Omega)\ne\mbox{\rm co}\big(\mathbb{S}\cap N(0;\Omega)\big).
\end{equation*}
Let us first check that $0\in{\rm bd}(\Omega)$. Indeed, fix any $\ve>0$ and find $k\in\N$ with $1/k<\ve/2$. This implies that $(\ve/2)e_k\in\B(0;\ve)\cap\Omega^c$, and hence $0\in\mbox{\rm bd}(\Omega)$ since $\ve>0$ was chosen arbitrarily. Due to $\pm(1/k)e_k\in\Omega$ for all $k\in\N$, we see that the inclusion $x\in N(0;\Omega)$ yields $\la x,e_k\ra=0$ as $k\in\N$. This shows that $N(0;\Omega)=\{0\}$, and thus $\mathbb{S}\cap N(0;\Omega)=\emp$. On the other hand, we know that $\Hat d(\cdot;\Omega)$ is (Lipschitz) continuous and convex, and therefore $\partial\Hat d(0;\Omega)\ne\emp$ by basic convex analysis. This tells us that formula \eqref{SDF} fails to hold in Hilbert spaces.}
\end{Example}\vspace*{-0.3in}

\section{Concluding Remarks}\label{sec:concl}\vspace*{-0.1in}

This paper provides a unified study of the class of minimal time functions, including their convex generalized differentiation, in the framework of arbitrary LCTV spaces by using a new notion of the tangent set closure with respect to constant dynamics. Our further attention is paid to the class of signed minimal time functions, which has not been considered in the literature, while being certainly important for applications. We establish general properties of the signed minimal time functions in LCTV spaces and then propose a variational approach to their generalized differentiation concentrating on a special class of convex signed distance functions in finite dimensions. Our variational approach and obtained results open the door to future investigation and applications of the new class of signed minimal time functions in convex and nonconvex settings of both finite- and infinite-dimensional spaces.
\end{document}